\documentclass{amsart}

\usepackage{amssymb}
\usepackage{array}

\allowdisplaybreaks

\newtheorem{theorem}{Theorem}[section]
\newtheorem{lemma}[theorem]{Lemma}

\theoremstyle{definition}
\newtheorem{definition}[theorem]{Definition}
\newtheorem{remark}[theorem]{Remark}

\newcommand{\Z}{\mathbb{Z}}
\newcommand{\R}{\mathbb{R}}
\newcommand{\C}{\mathbb{C}}

\newcommand{\Pol}{\mathbb{P}}

\newcommand{\Po}{\mathrm{P}}

\newcommand{\Di}{\mathrm{D}}

\begin{document}

\title[Taylor series for the Askey-Wilson operator]
{Taylor series for the Askey-Wilson operator \\ and classical
summation formulas}

\author[B. L\'opez, J.M. Marco and J. Parcet]
{Bernardo L\'opez, Jos\'{e} Manuel Marco and Javier Parcet}

\address{Department of Mathematics, Universidad Aut\'{o}noma de
Madrid.}

\email{bernardo.lopez@uam.es} \email{javier.parcet@uam.es}


\begin{abstract}
An analogue of Taylor's formula, which arises by substituting the
classical derivative by a divided difference operator of
Askey-Wilson type, is developed here. We study the convergence of
the associated Taylor series. Our results complement a recent work
by Ismail and Stanton. Quite surprisingly, in some cases the
Taylor polynomials converge to a function which differs from the
original one. We provide explicit expressions for the integral
remainder. As application, we obtain some summation formulas for
basic hypergeometric series. As far as we know, one of them is
new. We conclude by studying the different forms of the binomial
theorem in this context.
\end{abstract}

\maketitle

\section{Introduction and definitions}
\label{Section1}

The problem of expanding a function with respect to a given
polynomial basis has many implications in analysis. The simplest
example of this kind is the Taylor's expansion theorem. In this
paper, we replace the classical derivative by a difference
operator of Askey-Wilson type. Our results complement the paper
\cite{IS2} of Ismail and Stanton and are a natural continuation of
the point of view presented in \cite{MP}, where a new approach to
the theory of classical hypergeometric polynomials is given. In
contrast with \cite{IS2}, our aim is to find sufficient conditions
for the Taylor series to converge, but not necessarily to the
original function. In this more general setting, we may consider
non-necessarily entire functions and we give an explicit
expression for the limit of the remainders in terms of a contour
integral. Using this and a new estimate for the $q$-shifted
factorials, which might be of independent interest, we obtain a
summation formula which is new as far as we know. As we explain
below, it can be regarded as a \emph{non-symmetrized} version of
the non-terminating $q$-Saalsch\"{u}tz sum. As applications, we
also provide a new proof of the $q$-Gauss summation formula and a
list of binomial type summation formulas in the same line than
Ismail's paper \cite{I}.

\vskip3pt

Now we give some definitions which will be used in what follows.
The notions we are presenting were already introduced in \cite{MP}
with the aim of studying some aspects of the theory of
hypergeometric polynomials. The relevance of this approach is
justified in \cite{MP}, where a more detailed exposition is given.

\begin{definition}
Let $\mathcal{P}$ be the set of quadratic symmetric polynomials
$\Po(x,y)$, normalized so that $\Po(x,y) = x^2 + y^2 - 2axy -
2b(x+y) + c$ with $a,b,c$ complex numbers. A sequence of complex
numbers $(x_t)$, with the index $t$ running over the set
$\frac{1}{2} \Z = \{k/2: \, k \in \Z\}$, will be called a
$\Po$-\textbf{sequence} if $$\textstyle \Po(x_t,y) = (y - x_{t +
\frac{1}{2}}) (y - x_{t - \frac{1}{2}}) \qquad \mbox{for all}
\qquad t \in \frac{1}{2} \Z.$$
\end{definition}

Recalling the natural action of the affine group of the complex
plane $\mbox{Aff}(\C)$ on the set $\mathcal{P}$, we can consider
the orbits of this action. As it becomes clear in \cite{MP}, these
orbits provide a very natural classification in the theory. In
particular, each symmetric polynomial in $\mathcal{P}$ can be
rewritten, under affine transformations, in a \textbf{canonical
form}. In the following table, we summarize the canonical forms
and the corresponding $\Po$-sequences which we shall use in our
study. We shall also use the parameter $\lambda \neq 0$ defined by
the relation $a = \frac{1}{2} (\lambda + \lambda^{-1})$. Here we
recall that the basis $q$ will be given by $\lambda^2$.
\begin{center}
\setlength{\extrarowheight}{4pt}
\begin{tabular}{|r|c|c|} \multicolumn{1}{c}{} &
\multicolumn{1}{c}{\textbf{Canonical form}} &
\multicolumn{1}{c}{\textbf{$\Po$-sequence}} \\ \cline{1-3}
\textbf{T} & $x^2 + y^2 - 2 a xy + a^2 - 1$ & $\frac{1}{2}
(\lambda^{2t} u + \lambda^{-2t} u^{-1})$ \\ \cline{1-3} \textbf{G}
& $x^2 + y^2 - 2 a xy$ & $\lambda^{2t} u$ \\ \cline{1-3}
\textbf{Q} & $(x-y)^2 - \frac{1}{2}(x+y) + \frac{1}{16}$ & $t^2 +
2tu + u^2$ \\ \cline{1-3} \textbf{A} & $(x-y)^2 - \frac{1}{4}$ &
$t + u$ \\ \cline{1-3} \textbf{C} & $(x-y)^2$ & $u$
\\ \cline{1-3}
\end{tabular} \\ \vskip10pt \textsc{Table I.} Canonical
$\Po$-sequences. \vskip6pt
\end{center}

In the first two rows we assume $a \neq \pm 1$. The capital
letters in the left column are acronyms of the names we adopted in
\cite{MP} for the canonical forms: \textbf{Trigonometric},
\textbf{Geometric}, \textbf{Quadratic}, \textbf{Arithmetic} and
\textbf{Continuous}. We shall say that $x_0$ is the base point of
the $\Po$-sequence $(x_t)$. Obviously, $\Po(x,y) = y^2 - 2A(x)y +
B(x)$, where $A(x) = ax + b$ and $B(x) = x^2 - 2bx + c$. This
allows us to consider the discriminant of $\Po$ as a function of
the variable $x$ which, up to a constant factor, is given by
$$\delta(x) = (a^2-1)x^2 + 2b(a+1)x + b^2 - c.$$

\begin{remark}
$\Po$-sequences arise from the recurrence $x_{t \pm \frac{1}{2}} =
A(x_t) \pm \sqrt{\delta(x_t)}$. In particular, given a complex
number $\xi$, there are at most two $\Po$-sequences with base
point $\xi$, one for each choice of the sign for the square root
$\sqrt{\delta(\xi)}$.
\end{remark}

\begin{definition} \label{Definicion-Phik}
For any complex number $x$, we consider a $\Po$-sequence $(x_t)$
with base point $x$. Then we define $\Phi_0(x,y) = 1$ and the
\textbf{polynomials} $\Phi_k(x,y)$, for any positive integer $k
\ge 1$, as follows $$\Phi_k(x,y) = \prod_{j=0}^{k-1} (y - x_{j -
\frac{k-1}{2}}).$$
\end{definition}

\begin{remark}
Recall that $\Phi_k(x,\cdot)$ does not depend on the chosen
$\Po$-sequence $(x_t)$.
\end{remark}

\section{Taylor series}
\label{Section2}

Given an open subset $\Omega$ of the complex plane, we denote the
space of analytic functions in $\Omega$ by $\mathcal{H}(\Omega)$.
Also, by $\gamma \simeq 0 \ (\mbox{mod} \, \Omega)$ we mean that
$\gamma$ is a cycle in $\Omega$ homologous to zero with respect to
$\Omega$. Finally, given $z \in \Omega$, $\mbox{Ind}(\gamma,z)$
denotes the index of $z$ with respect to $\gamma$.

\begin{lemma} \label{Lema-Cauchy}
Given an open subset $\Omega$ of the complex plane, let us denote
by $\Delta_{\Omega}$ the diagonal of $\Omega \times \Omega$. Then,
for any $f \in \mathcal{H}(\Omega)$, the function $$f_d(u,v) =
\frac{f(u)-f(v)}{u-v} \qquad \mbox{for} \qquad (u,v) \in (\Omega
\times \Omega) \setminus \Delta_{\Omega}$$ can be continuously
extended to an analytic function $f_d: \Omega \times \Omega
\rightarrow \C$. Moreover, if $\gamma \simeq 0 \ (\textnormal{mod}
\, \Omega)$ and $\textnormal{Ind}(\gamma,u) =
\textnormal{Ind}(\gamma,v) = 1$, then we have $$f_d(u,v) =
\frac{1}{2 \pi i} \int_{\gamma} \frac{f(y)}{(y-u)(y-v)} \, dy.$$
\end{lemma}

\begin{proof}
It is a simple consequence of Cauchy's integral formula.
\end{proof}

Let us consider the set $\Omega_1 = \{x \in \C: \, A(x) \pm
\sqrt{\delta(x)} \in \Omega\}$. By Lemma \ref{Lema-Cauchy}, given
a quadratic symmetric polynomial $\Po \in \mathcal{P}$, we define
the \textbf{Askey-Wilson operator} $\Di: \mathcal{H}(\Omega)
\rightarrow \mathcal{H}(\Omega_1)$ as follows
\begin{equation} \label{Ecuacion-D-Analitico}
\Di f(x) = f_d \big( A(x) + \sqrt{\delta(x)}, A(x) -
\sqrt{\delta(x)} \big) = \frac{1}{2 \pi i} \int_{\gamma}
\frac{f(y)}{\Po(x,y)} \, dy,
\end{equation}
with $\gamma \simeq 0 \ (\mbox{mod} \, \Omega)$ and
$\mbox{Ind}(\gamma, A(x) \pm \sqrt{\delta(x)}) = 1$. On the other
hand, if we define recursively the sets $\Omega_k$ by
$\Omega_{k+1} = \{ x \in \C : \, A(x) \pm \sqrt{\delta(x)} \in
\Omega_k \}$ with $\Omega_0 = \Omega$, we may consider the
iterated operators $$\Di^k: \mathcal{H}(\Omega) \to
\mathcal{H}(\Omega_k).$$

\begin{lemma} \label{Lema-Dk-Analitico}
Let $\Omega$ be an open set of the complex plane and $\gamma
\simeq 0 \ (\textnormal{mod} \, \Omega)$. Then, given $f \in
\mathcal{H}(\Omega)$ and $x \in \Omega_k$, we have the following
identity $($with the obvious limits for $\lambda = \pm 1)$ $$\Di^k
f (x) = \Big( \prod_{j=0}^{k-1} \frac{\lambda^{k-j} -
\lambda^{j-k}}{\lambda - \lambda^{-1}} \Big) \, \frac{1}{2 \pi i}
\int_{\gamma} \frac{f(y)}{\Phi_{k+1}(x,y)} \, dy \qquad \mbox{for}
\qquad k \ge 0$$ if $x$ is the base point of a $\Po$-sequence
$(x_t)$ with $\textnormal{Ind}(\gamma, x_{j - \frac{k}{2}}) = 1$
for $j = 0,1, \ldots, k$. In particular, we shall consider the
operator $$\partial_k: \mathcal{H}(\Omega) \rightarrow
\mathcal{H}(\Omega_k) \quad \mbox{defined by} \quad \partial_k f
(x) = \frac{1}{2 \pi i} \int_{\gamma} \frac{f(y)}{\Phi_{k+1}(x,y)}
\, dy.$$
\end{lemma}

\begin{proof}
The cases $k=0$ and $k=1$ follow from Cauchy's integral formula
and (\ref{Ecuacion-D-Analitico}) respectively. Therefore, the
general case follows by induction from the following relation
$$\Di \Big[ \frac{1}{\Phi_k(\cdot,y)} \Big] (x) =
\frac{\Phi_k(x_{1/2},y)^{-1} - \Phi_k(x_{-1/2},y)^{-1}}{x_{1/2} -
x_{-1/2}} = \frac{\lambda^k - \lambda^{-k}}{\lambda -
\lambda^{-1}} \Phi_{k+1}(x,y)^{-1},$$ which is not difficult to
check with the aid of the identity $$\frac{x_{k/2} -
x_{-k/2}}{x_{1/2} - x_{-1/2}} = \frac{\lambda^k -
\lambda^{-k}}{\lambda - \lambda^{-1}}.$$ This identity was proved
in \cite{MP} for any $\Po$-sequence. This completes the proof.
\end{proof}

The following result provides the sequence of Taylor polynomials
associated to a given function $f \in \mathcal{H}(\Omega)$ and the
corresponding remainder term, with respect to the Askey-Wilson
operator.

\begin{theorem} \label{Teorema-Taylor-Series}
Let $\Omega$ be an open subset of the complex plane and $\gamma
\simeq 0 \ (\textnormal{mod} \, \Omega)$. Then, given an analytic
function $f \in \mathcal{H}(\Omega)$ and $x \in \Omega$ with
$\textnormal{Ind}(\gamma, x) = 1$, we can recover $f(x)$ as
$$f(x) = \sum_{k=0}^n \partial_k f(z_{k/2}) \prod_{j=0}^{k-1} (x -
z_j) \, + \, \mathcal{R}_n f(x)$$ where $\mathcal{R}_n f(x)$ is
given by $$\mathcal{R}_n f(x) = \frac{1}{2 \pi i} \int_{\gamma}
\frac{f(y)}{y-x} \prod_{j=0}^n \frac{x-z_j}{y-z_j} \, dy,$$ with
$(z_t)$ a $\Po$-sequence such that $z_j \in \Omega$ and
$\textnormal{Ind}(\gamma,z_j) = 1$ for $j = 0,1, \ldots, n$.
\end{theorem}

\begin{proof}
The relation below follows by induction on $n$ $$\frac{1}{y-x} =
\frac{\Phi_{n+1}(z_{n/2},x)}{(y-x) \Phi_{n+1}(z_{n/2},y)} +
\sum_{k=0}^n \frac{\Phi_k (z_{(k-1)/2},x)}{\Phi_{k+1}
(z_{k/2},y)}.$$ Multiplying by $(2 \pi i)^{-1} f(y)$, we are done
by Lemma \ref{Lema-Dk-Analitico} and Cauchy's formula.
\end{proof}

\begin{remark} \label{Observacion-Taylor-Clasica}
When dealing with the continuous form \textbf{C}, it turns out
that the Askey-Wilson operator can be regarded as the classical
derivative, see \cite{MP} for the details. In particular, Theorem
\ref{Teorema-Taylor-Series} reduces to the classical Taylor series
$$f(x) = \sum_{k=0}^n \frac{f^{(k)}(z)}{k!} (x-z)^k + f_{n+1}(x)
(x-z)^{n+1},$$ where $$f_{n+1}(x) = \frac{1}{2 \pi i}
\int_{\gamma} \frac{f(y)}{(y-z)^{n+1}(y-x)} \, dy.$$
\end{remark}

\begin{remark} \label{Observacion-Polinomio-Interpolador}
If the points $z_j$ are pairwise distinct for $j = 0,1, \ldots,
n$, we can use Lemma \ref{Lema-Dk-Analitico} to express
$\partial_k f (z_{k/2})$ as a sum of residues and Theorem
\ref{Teorema-Taylor-Series} gives $$f(x) - \mathcal{R}_n f(x) =
\sum_{k=0}^n \Big( \sum_{j=0}^k \frac{f(z_j)}{\prod_{i \neq j}
(z_j - z_i)} \Big) \prod_{j=0}^{k-1} (x - z_j),$$ which is
Newton's divided difference formula for the interpolation
polynomial. In particular, Theorem \ref{Teorema-Taylor-Series}
holds for any collection $z_0, z_1, \ldots, z_n$ of pairwise
distinct points. However, as we shall see below, the relevance of
Theorem \ref{Teorema-Taylor-Series} lies in the established
connection with the Askey-Wilson operator.
\end{remark}

Now we study the convergence of the Taylor series for
$\Po$-sequences $(z_t)$ with $z_0, z_1, z_2, \ldots$ bounded. This
covers the geometric canonical form for $|q| < 1$ and the
continuous form, whose associated $\Po$-sequences are constant.
Given $\mathrm{r} > 0$, we shall denote by
$\mathbb{D}_{\mathrm{r}} = \{z \in \C: \, |z| < \mathrm{r}\}$ the
open disk of radius $\mathrm{r}$. As it was pointed out by the
referee, the result below is closely related to Wallisser's paper
\cite{W}.

\begin{theorem} \label{Teorema-Convergencia-Series}
Let $(z_t)$ be a $\Po$-sequence satisfying that the subsequence
$z_0, z_1, z_2,..$ is bounded. Then, given $\mathbf{z} =
\limsup_{k \ge 0} |z_k|$, $\delta > 0$ and $f \in
\mathcal{H}(\mathbb{D}_{\mathbf{r} + \delta})$ with $\mathbf{r} >
2 \mathbf{z}$, the Taylor polynomials $$\mathcal{S}_n f(x) =
\sum_{k=0}^n \partial_k f(z_{k/2}) \prod_{j=0}^{k-1} (x - z_j)$$
converge uniformly to $f(x)$ as $n \rightarrow \infty$ on the disk
$\mathbb{D}_{\mathbf{x}}$ for any $\mathbf{x} < \mathbf{r} - 2
\mathbf{z}$.
\end{theorem}

\begin{proof}
By Theorem \ref{Teorema-Taylor-Series}, it suffices to check
$\mathcal{R}_n f \rightarrow 0$ uniformly on the disk
$\mathbb{D}_{\mathbf{x}}$. Let us take $\gamma$ to be a
circumference centered at $0$ with radius $\mathbf{r}
> \mathbf{x} + 2 \mathbf{z}$. It turns out that, if
$\mathrm{M}_{\mathbf{r}}(f)$ denotes the supremum of $|f(y)|$ when
$y$ runs over $\gamma$, we have $$|\mathcal{R}_n f(x)| \le
\mathbf{r} \frac{\mathrm{M}_{\mathbf{r}}(f)}{\mathbf{r} -
\mathbf{x}} \ \Big( \frac{\mathbf{x} + \mathbf{z}}{\mathbf{r} -
\mathbf{z}} \Big)^{n+1} \longrightarrow 0 \qquad \mbox{as} \qquad
n \rightarrow \infty.$$ Since the given bound holds for every $x
\in \mathbb{D}_\mathbf{x}$, the proof is complete.
\end{proof}

The following is a convergence theorem with integral remainder.
This case is applicable to the trigonometric and quadratic
canonical forms as well as to the geometric canonical form with
$|q| > 1$. Let $(z_t)$ be a $\Po$-sequence satisfying $\sum_{k \ge
0} 1/|z_k| < \infty$. Then, the product
$$\mathrm{H}(x) = \prod_{j=0}^\infty \big( 1 - x/z_j \big)$$ is
an entire function of $x$ with zeros at $x=z_j$ for $j=0,1,2,
\ldots$

\begin{theorem} \label{Theorem-Integral-Remainder}
Let $(z_t)$ be a $\Po$-sequence satisfying that the subsequence
$z_0, z_1, z_2, \ldots$ lies in $\R_- = \{ x \in \R: x < 0 \}$ and
$\sum_{k \ge 0} 1/|z_k| < \infty$. Assume that $f$ is analytic in
an open neighborhood of the left half-plane $$\Omega = \Big\{ z
\in \C : \, \mathrm{Re} \, z \le 0 \Big\},$$ and is controlled by
\begin{equation} \label{Control}
|f(z)| \le \mathrm{C} \big( 1 + |z| \big)^{\mathrm{M}} \quad
\mbox{for all} \quad z \in \Omega
\end{equation}
and some absolute constants $\mathrm{C}$ and $\mathrm{M}$. In that
case the pointwise limit
$$\mathcal{S}_{\infty}f(x) = \lim_{n \to \infty} \mathcal{S}_n
f(x)$$ exists when $\mathrm{Re} \, x < 0$ and we have
$$f(x) = \mathcal{S}_{\infty} f(x) + \mathcal{R}_\infty f(x)$$
with $$\mathcal{R}_\infty f(x) = \frac{1}{2 \pi i} \int_{-i
\infty}^{i \infty} \frac{f(y)}{y-x}
\frac{\mathrm{H}(x)}{\mathrm{H}(y)} dy = \frac{1}{2 \pi i}
\int_{-i \infty}^{i \infty} \frac{f(y)}{y-x} \prod_{j=0}^{\infty}
\frac{x - z_j}{y - z_j} dy.$$
\end{theorem}

\begin{proof}
According to Theorem \ref{Teorema-Taylor-Series} we have
$$\mathcal{R}_n f(x) = \frac{1}{2 \pi i} \Big( \int_{-ir}^{ir}
\frac{f(y)}{y-x} \prod_{j=0}^n \frac{x-z_j}{y-z_j} \, dy +
\int_{\gamma_r} \frac{f(y)}{y-x} \prod_{j=0}^n \frac{x-z_j}{y-z_j}
\, dy \Big),$$ whenever $\mbox{Re} \, x < 0$ and where $r$ and
$\gamma_r: [0, \pi] \to \C$ are determined by
$$\max \Big\{|x|, |z_0|, |z_1|, \ldots, |z_n| \Big\} < r \quad
\mbox{and} \quad \gamma_r(t) = i r e^{it} \ (0 \le t \le \pi).$$
On the other hand, it easily follows from (\ref{Control}) that for
any $y \in \gamma_r[0,\pi]$ we have $$\Big| \frac{f(y)}{y-x}
\prod_{j=0}^n \frac{x-z_j}{y-z_j} \Big| \le \mathrm{C}' \big( 1 +
r \big)^{\mathrm{M} - n - 2},$$ for some constant $\mathrm{C}'$
depending on $x, z_0, z_1, \ldots, z_n$ but not on $y$. Therefore,
since we may take the parameter $r$ arbitrary large, we deduce
that for any integer $n \ge \mathrm{M}$ we have
$$\lim_{r \to \infty} \int_{\gamma_r} \frac{f(y)}{y-x} \prod_{j=0}^n
\frac{x-z_j}{y-z_j} \, dy = 0.$$ In other words, for $n \ge
\mathrm{M}$ we have $$\mathcal{R}_n f(x) = \frac{1}{2 \pi i}
\prod_{j=0}^n \big( 1 - x/z_j \big) \int_{-i \infty}^{i \infty}
\frac{f(y)}{y-x} \prod_{j=0}^n \big( 1 - y/z_j \big)^{-1} \, dy.$$
Then, estimating as above we find for $n \ge \mathrm{M}$ $$\Big|
\frac{f(y)}{y-x} \prod_{j=0}^n \big( 1 - y/z_j \big)^{-1} \Big|
\le \mathrm{C}'' \big( 1 + |y| \big)^{\mathrm{M}-1} \prod_{j=0}^n
\Big( 1 + \frac{|y|^2}{|z_j|^2} \Big)^{-1/2} \le \mathrm{C}'''
\big( 1 + |y| \big)^{-2}.$$ Therefore, the assertion follows from
the dominated convergence theorem.
\end{proof}

\begin{remark}
Similar arguments show that the convergence in Theorem
\ref{Theorem-Integral-Remainder} is uniform over the compact sets
of $\{x \in \C : \, \mbox{Re} \, x < 0 \}$. On the other hand, it
is clear that Theorem \ref{Theorem-Integral-Remainder} can be
restated by taking $z_0, z_1, z_2, \ldots$ in $\R_+$, $$\Omega =
\Big\{ z \in \C : \, \mbox{Re} \, z \ge 0 \Big\} \quad \mbox{and}
\quad \mathcal{R}_\infty f(x) = \frac{-1}{2 \pi i} \int_{-i
\infty}^{i \infty} \frac{f(y)}{y-x} \prod_{j=0}^{\infty} \frac{x -
z_j}{y - z_j} \, dy.$$ This shows that the quadratic
$\Po$-sequences are contemplated by Theorem
\ref{Theorem-Integral-Remainder}. We have emphasized this point
following a suggestion from the referee. However, when dealing
with the trigonometric form, we might consider a growth
restriction weaker than (\ref{Control}). We omit the proof for
lack of space and refer the reader to \cite{IS2}.
\end{remark}

\begin{remark}
The main result in Ismail/Stanton's paper \cite[Theorem 3.1]{IS2}
considers Taylor series of entire functions with respect to the
trigonometric canonical form. Our main contribution in Theorem
\ref{Theorem-Integral-Remainder} is that we do not require $f$ to
be entire. Note that the regularity and growth restrictions in
\cite{IS2} allowed the authors to obtain the Taylor series of
\emph{products} of $q$-shifted factorials. Theorem
\ref{Theorem-Integral-Remainder} and Lemma \ref{Lemma-Bernardo}
below allow us to consider in the next section \emph{quotients} of
$q$-shifted factorials.
\end{remark}

\section{Summation formulas}
\label{Section3}

Now we compute the Taylor coefficients of a couple of functions
made up of $q$-shifted factorials. Then we analyze the convergence
of the corresponding Taylor series. We obtain a new proof of the
$q$-Gauss summation formula and, as far as we know, a new
summation formula.

\subsection{Geometric case}
\label{Subsection-Geometric}

We begin by applying Theorem \ref{Teorema-Convergencia-Series} to
a particular case. This will provide a new proof of the $q$-Gauss
summation formula. We shall work with the Askey-Wilson operator
which arises from the geometric canonical form
$$\Di f(z) = \frac{f(q^{1/2} z) - f(q^{-1/2} z)}{q^{1/2}z -
q^{-1/2}z}.$$ Given $\alpha, \beta \in \C$ and $|q| < 1$, we
consider the function $$f(x) = \frac{(\alpha x;q)_{\infty}}{(\beta
x;q)_{\infty}} = \prod_{j=0}^{\infty} \frac{1 - q^j \alpha x}{1 -
q^j \beta x}.$$ By induction, it is not difficult to check that
$$\partial_k f(z) = \Big( \prod_{j=0}^{k-1} \frac{\beta - q^j
\alpha}{1 - q^{j+1}} \Big) \frac{(\alpha z
q^{k/2};q)_{\infty}}{(\beta z q^{-k/2};q)_{\infty}}.$$ Therefore,
taking $\xi \in \C$ and $z_t = q^t \xi$ in Theorem
\ref{Teorema-Convergencia-Series}, we obtain
$$\frac{(\alpha x;q)_{\infty} (\beta \xi;q)_{\infty}}{(\beta
x;q)_{\infty} (\alpha \xi;q)_{\infty}} = \sum_{k=0}^{\infty}
\frac{(\xi/x;q)_k}{(q;q)_k(\alpha \xi;q)_k} \Big(
\prod_{j=0}^{k-1} (\beta - q^j \alpha) \Big) x^k.$$ In particular,
with the usual notation for basic hypergeometric series
$$_2\phi_1 \Big[ \begin{array}{c} \xi/x, \ \alpha / \beta \\
\alpha \xi \end{array} ; \ q, \ \beta x \Big] = \frac{(\alpha x,
\beta \xi ;q)_{\infty}}{(\beta x, \alpha \xi ;q)_{\infty}}.$$ In
other words, we have obtained the $q$-Gauss classical summation
formula.

\subsection{An estimate for $q$-shifted factorials}

In this paragraph we provide an estimate for $q$-shifted
factorials that will be needed below. Although we just need a
weaker estimate, we have included the sharpest result we know
since it might be of independent interest.

\begin{lemma} \label{Lemma-Bernardo}
The following assertions hold for $0 < q < 1$:
\begin{itemize}
\item[(a)] There exists some $\mathrm{C}_q > 0$ such that every $x
\in \C$ satisfies $$\Big| \prod_{j=0}^\infty \big( 1 - q^j x \big)
\Big| \le \mathrm{C}_q |x|^{\frac12 + \frac{\log |x|}{2 \log
q^{-1}}}.$$

\item[(b)] Let $\mathcal{A}$ be any closed set in $\C$ with $q
\mathcal{A} \subset \mathcal{A}$ and $q^{-j} \notin \mathcal{A}$
for $j = 0,1,2, \ldots$ Then, given any $\delta > 0$ there exists
a constant $\mathrm{C}_q(\mathcal{A}, \delta) > 0$ such that for
all $x \in \mathcal{A} \setminus \mathbb{D}_\delta$, we have
$$\Big| \prod_{j=0}^\infty \big( 1 - q^j x \big) \Big| \ge
\mathrm{C}_q(\mathcal{A},\delta) |x|^{\frac12 + \frac{\log |x|}{2
\log q^{-1}}}.$$
\end{itemize}
\end{lemma}

\begin{proof}
Let us write $$h(x) = \prod_{j=0}^\infty \big( 1 - q^j x \big)
\quad \mbox{and} \quad h_n(x) = \prod_{j=0}^{n-1} \big( 1 - q^j x
\big).$$ We begin by proving the upper estimate. Since the
function $|x|^{\frac12 + \frac{\log |x|}{2 \log q^{-1}}}$ explotes
at $x=0$, we may assume that $|x| > 1$. Then we fix the only
integer $m$ satisfying the estimates $q^m |x| \le 1 < q^{m-1} |x|$
and observe that
\begin{equation} \label{Star}
\big| h(x) \big| = \Big| h(q^m x) \prod_{j=0}^{m-1} \big( 1 - q^j
x \big) \Big| = |x|^m q^{{{m} \choose {2}}} \big| h(a) h_m(b)
\big|,
\end{equation}
where $a = q^m x$ and $b = q^{-m+1} x^{-1}$. Our choice of $m$
implies $$m = \frac{\log |x|}{\log q^{-1}} + \rho \quad \mbox{for
some} \quad 0 \le \rho < 1.$$ Moreover, a simple computation gives
\begin{equation} \label{Star2}
|x|^m q^{{{m} \choose {2}}} = q^{\frac{\rho(\rho-1)}{2}}
|x|^{\frac{1}{2} + \frac{\log |x|}{2 \log q^{-1}}}.
\end{equation}
Therefore, since $|a| \le 1$ and $q \le |b| < 1$, we have $|h(a)
h_n(b)| \le h(-1)^2$ and we deduce the assertion with constant
$\mathrm{C}_q = q^{-1/8} h(-1)^2$. Now we prove the lower
estimate. Let us consider the set $$\mathcal{B} = \Big\{ z \in \C:
\, z^{-1} \in \mathcal{A}, \, q \le |z| \le 1 \Big\}.$$ The
function $h$ does not vanish on $\mathcal{A} \cup \mathcal{B}$.
Therefore, since $h_n \to h$ uniformly over compact sets, we
deduce that $\min \big\{ |h_n(z)|, |h(z)| \big\} \ge
\mathrm{C}_q(\mathcal{A})$ for some $\mathrm{C}_q(\mathcal{A})$,
for all $n \ge 1$ and all $z \in \mathcal{A} \cup \mathcal{B}$
with $|z| \le 1$. Thus, if $|x|
> 1$ we deduce from (\ref{Star}) and (\ref{Star2}) (recall that $a
\in \mathcal{A}$ and $b \in \mathcal{B}$) that
$$\Big| \prod_{j=0}^\infty \big( 1 - q^j x \big) \Big| \ge
\mathrm{C}_q(\mathcal{A}) |x|^{\frac12 + \frac{\log |x|}{2 \log
q^{-1}}} \quad \mbox{for all} \quad x \in \mathcal{A} \setminus
\overline{\mathbb{D}}_1.$$ The assertion of $x \in \mathcal{A}$
with $\delta \le |x| \le 1$ follows by replacing
$\mathrm{C}_q(\mathcal{A})$ by $\mathrm{C}_q(\mathcal{A},
\delta)$.
\end{proof}

\begin{remark}
Let us observe how to construct sets of type $\mathcal{A}$. Since
$q^{-j} \notin \mathcal{A}$ for $j = 0,1,2, \ldots$ and
$\mathcal{A}$ is closed, there must exists an open neighborhood
$\mathcal{U}$ of $1$ such that $\mathcal{A} \cap \mathcal{U} =
\emptyset$. Now, recalling that $q \mathcal{A} \subset
\mathcal{A}$ if and only if $q^{-1} \mathcal{A}^c \subset
\mathcal{A}^c$, we must have that the neighborhood of $q^{-j}$
defined by $q^{-j} \mathcal{U}$ is contained in $\mathcal{A}^c$.
Thus, the \emph{biggest} $\mathcal{A}$'s are sets generated by a
neighborhood $\mathcal{U}$ of $1$ in the following way
$$\mathcal{A} = \bigcap_{j=0}^\infty \big( q^{-j} \mathcal{U}
\big)^c.$$
\end{remark}

\subsection{Trigonometric case}
\label{Subsection-Trigonometric}

In this paragraph we apply Theorem
\ref{Theorem-Integral-Remainder} to $$f_{\alpha \beta} (x) =
f_{\alpha \beta} \Big( \frac{u + u^{-1}}{2} \Big) = \frac{(\alpha
u; q)_\infty (\alpha/u ;q)_{\infty}}{(\beta u;q)_\infty
(\beta/u;q)_{\infty}} \quad \mbox{with} \quad x = \frac12 (u +
u^{-1})$$ and $(\alpha, \beta) \in \R \times \R_+$. As a function
of $x$, $f_{\alpha,\beta}$ is well-defined since $(\gamma u,\gamma
/ u;q)_{\infty}$ is symmetric under the mapping $u \mapsto
u^{-1}$. Moreover, it is analytic in a neighborhood of the left
half-plane (as a function of $x$) since $(\gamma u,\gamma /
u;q)_{\infty}$ is analytic in $\C \setminus \{0\}$ for any $\gamma
\in \C$ and has zeros in $\R_+$ when $\gamma \in \R_+$. We shall
use the Askey-Wilson operator which arises from the trigonometric
canonical form $$\Di f \Big( \frac{u + u^{-1}}{2} \Big) =
\frac{2}{(\lambda - \lambda^{-1})(u - u^{-1})} \left[ f \Big(
\frac{\lambda u + \lambda^{-1} u^{-1}}{2} \Big) - f \Big(
\frac{\lambda^{-1} u + \lambda u^{-1}}{2} \Big) \right]$$ with
$q=\lambda^2$ and $0 < \lambda < 1$. Grouping the common factors,
it can be checked that
$$\Di f_{\alpha \beta} \Big( \frac{u + u^{-1}}{2} \Big) = 2 \,
\frac{\alpha - \beta}{q-1} \ f_{\lambda \alpha, \lambda^{-1}
\beta} \Big( \frac{u + u^{-1}}{2} \Big).$$ Moreover, by induction
on $k$ we obtain $$\partial_k f_{\alpha \beta} \Big( \frac{u +
u^{-1}}{2} \Big) = \frac{2^k}{(q;q)_k} \ \prod_{j=0}^{k-1} (\beta
- q^j \alpha) \ f_{\lambda^k \alpha, \lambda^{-k} \beta} \Big(
\frac{u + u^{-1}}{2} \Big).$$ In particular, taking $z_t =
\frac{1}{2} \big( q^t \xi + q^{-t} \xi^{-1} \big)$ with $\xi < 0$,
we can write
$$\partial_k f_{\alpha \beta} (z_{k/2}) = \frac{(2 \beta)^k
(\alpha/\beta;q)_k}{(q,\alpha \xi,\beta/q^k \xi;q)_k} \ f_{\alpha
\beta} \Big( \frac{\xi + \xi^{-1}}{2} \Big).$$ According to
Theorem \ref{Teorema-Taylor-Series} and
\begin{equation} \label{Ecuacion-Factorizacion}
\frac{\zeta + \zeta^{-1}}{2} - \frac{\eta + \eta^{-1}}{2} =
\frac{\big( (\zeta \eta)^{1/2} - (\zeta \eta)^{-1/2} \big) \big(
(\zeta / \eta)^{1/2} - (\zeta / \eta)^{-1/2} \big)}{2},
\end{equation}
the Taylor polynomials of $f_{\alpha, \beta}$ are
$$\mathcal{S}_n f_{\alpha \beta} \Big( \frac{u + u^{-1}}{2} \Big)
= f_{\alpha \beta} \Big( \frac{\xi + \xi^{-1}}{2} \Big) \
\sum_{k=0}^n \frac{(\alpha / \beta,\xi u ,\xi / u ;q)_k}{(q,\alpha
\xi,q \xi / \beta;q)_k} \ q^k.$$ In other words, in terms of
Theorem \ref{Theorem-Integral-Remainder} we have
\begin{equation} \label{Sinfty}
\mathcal{S}_\infty f_{\alpha \beta} \Big( \frac{u + u^{-1}}{2}
\Big) = \frac{(\alpha \xi,\alpha/\xi ;q)_{\infty}}{(\beta
\xi,\beta/\xi;q)_{\infty}} \, _3\phi_2 \Big[
\begin{array}{c} \alpha/\beta, \ \xi u, \ \xi/u \\ \alpha \xi, \ q
\xi/ \beta \end{array}; \ q, \ q \Big].
\end{equation}
On the other hand, the subsequence $z_0, z_1, z_2, \ldots$ lies in
$\R_-$ with $\sum_{k \ge 0} 1 / |z_k| < \infty$. Thus, in order to
apply Theorem \ref{Theorem-Integral-Remainder}, it suffices to
check whether the bound (\ref{Control}) holds for some $\mathrm{M}
> 0$. To that aim we shall assume that $|u| \ge 1$ and $$|x| \ge 2
\max \Big\{ |\alpha|, |\alpha|^{-1}, \beta, \beta^{-1} \Big\}.$$
This implies $\frac14 |u| \le |x| \le |u|$ so that $|\alpha/u|,
\beta/|u| \le 1/2$ and the functions $(\alpha/u;q)_\infty$,
$(\beta/u;q)_\infty$ are bounded from above and below. Moreover,
according to Lemma \ref{Lemma-Bernardo}, we obtain the following
estimate $$\big| (\alpha u;q)_\infty \big| \le
\mathrm{C}_q(\alpha) |u|^{\frac12 + \frac{\log |u|}{2\log q^{-1}}}
|\alpha|^{\frac{\log |u|}{2\log q^{-1}}} |u|^{\frac{\log
|\alpha|}{2\log q^{-1}}} = \mathrm{C}_q(\alpha) |u|^{\frac12 +
\frac{\log|u|}{2 \log q^{-1}} + \frac{\log |\alpha|}{\log
q^{-1}}},$$ with $\mathrm{C}_q(\alpha) = \mathrm{C}_q
|\alpha|^{\frac12 + \frac{\log |\alpha|}{2\log q^{-1}}}$.
Similarly, when $\mbox{Re} \, u \le 0$ (equivalently $\mbox{Re} \,
x \le 0$) $$\big| (\beta u;q)_\infty \big| \ge \mathrm{C}_q(\beta)
|u|^{\frac12 + \frac{\log |u|}{2\log q^{-1}}} |\beta|^{\frac{\log
|u|}{2\log q^{-1}}} |u|^{\frac{\log |\beta|}{2\log q^{-1}}} =
\mathrm{C}_q(\beta) |u|^{\frac12 + \frac{\log|u|}{2 \log q^{-1}} +
\frac{\log |\beta|}{\log q^{-1}}}.$$ In summary, we have seen that
$$\big| f_{\alpha, \beta} (x) \big| \le \mathrm{C}_q(\alpha,\beta)
|u|^{\mathrm{M}} \le 4 \mathrm{C}_q(\alpha,\beta) \big( 1 + |x|
\big)^{\mathrm{M}} \quad \mbox{with} \quad \mathrm{M} = \frac{\log
|\alpha/\beta|}{\log q^{-1}},$$ whenever $\mbox{Re} \, x \le 0$
and $|x| \ge 2 \max(|\alpha|, |\alpha|^{-1}, \beta, \beta^{-1})$.
However, the last restriction on $x$ can be dropped by continuity
and we are in the hypotheses of Theorem
\ref{Theorem-Integral-Remainder}. This gives the identity
$f_{\alpha,\beta}(x) = \mathcal{S}_\infty f_{\alpha,\beta}(x) +
\mathcal{R}_\infty f_{\alpha,\beta}(x)$ for $\mbox{Re} \, x < 0$.
Letting $f_\gamma(x) = (\gamma u, \gamma/u ; q)_\infty$ (so that
$f_{\alpha,\beta} = f_{\alpha}/f_{\beta}$), we have
$$\mathrm{H}(x) = \prod_{j=0}^\infty \Big( 1 - \frac{u+u^{-1}}{q^j
\xi + q^{-j} \xi^{-1}} \Big) = \prod_{j=0}^\infty \frac{(1- q^j\xi
u)(1 - q^j \xi / u)}{1+ q^{2j} \xi^2} = (-\xi^2;q^2)_\infty^{-1}
f_\xi(x).$$ Therefore, we have
\begin{equation} \label{Remainder}
\mathcal{R}_\infty f_{\alpha,\beta}(x) = \frac{1}{2\pi i} \int_{-i
\infty}^{i \infty} \frac{f_{\alpha,\beta}(y)}{y-x}
\frac{\mathrm{H}(x)}{\mathrm{H}(y)} \, dy = \frac{1}{2\pi i}
\int_{-i \infty}^{i \infty} \frac{f_{\alpha}(y)
f_{\xi}(x)}{f_{\beta}(y) f_{\xi}(y)} \, \frac{dy}{y-x}.
\end{equation}
Taking $y = \frac12 \big( v+v^{-1} \big)$, it follows from
(\ref{Sinfty}) and (\ref{Remainder}) that
\begin{eqnarray}
\label{New} \frac{(\alpha u,\alpha/u ;q)_{\infty}}{(\beta
u,\beta/u;q)_{\infty}} & = & \frac{(\alpha \xi,\alpha/\xi
;q)_{\infty}}{(\beta \xi,\beta/\xi;q)_{\infty}} \, _3\phi_2 \Big[
\begin{array}{c} \alpha/\beta, \ \xi u, \ \xi/u \\ \alpha \xi, \ q
\xi/ \beta \end{array}; \ q, \ q \Big] \\ \nonumber & + &
\frac{1}{2\pi i} \int_{-i \infty}^{i \infty} \frac{(\alpha v,
\alpha/v,\xi u, \xi /u; q)_\infty}{(\beta v, \beta/v,\xi v, \xi
/v; q)_\infty} \, \frac{(1-v^2)u}{(u-v)(uv-1)v} \, dv,
\end{eqnarray}
for $(\alpha, \beta) \in \R \times \R_+$, $\xi <0$ and $\mbox{Re}
\, u < 0$. This summation formula is new as far as we know. It can
be regarded as a \emph{non-symmetrized} version of the
non-terminating $q$-Saalsch\"{u}tz sum. To see this, we observe
that an obvious reformulation of Theorem
\ref{Theorem-Integral-Remainder} gives (\ref{New}) when $\mbox{Re}
\, x =0$ with $\mathcal{R}_\infty f_{\alpha,\beta}(x)$ replaced by
\begin{eqnarray*}
\mathcal{R}_\infty(x;\alpha,\beta,\xi) & = & \frac{1}{2 \pi i}
\Big( \int_{- i \infty}^{x-\delta i} + \int_{\gamma_{x}^+} +
\int_{x+\delta i}^{i \infty} \Big) \frac{f_\alpha(y)
f_\xi(x)}{f_\beta(y) f_\xi(y)} \, \frac{dy}{y-x},
\end{eqnarray*}
with $\gamma_x^+(t) = x - i \delta e^{it}$ for $0 \le t \le \pi$.
Transposing $(\beta,\xi)$ and changing signs
\begin{eqnarray*}
\mathcal{R}_\infty(-x;-\alpha,-\xi, -\beta) & = & \frac{1}{2 \pi
i} \Big( \int_{- i \infty}^{-x-\delta i} + \int_{\gamma_{-x}^+} +
\int_{-x+\delta i}^{i \infty} \Big) \frac{f_\alpha(-y)
f_\beta(x)}{f_\xi(-y) f_\beta(-y)} \, \frac{dy}{y+x} \\ & = &
\frac{1}{2 \pi i} \Big( \int_{i \infty}^{x+\delta i} +
\int_{\gamma_{x}^-} + \int_{x-\delta i}^{-i \infty} \Big)
\frac{f_\alpha(y) f_\beta(x)}{f_\beta(y) f_\xi(y)} \,
\frac{dy}{y-x},
\end{eqnarray*}
where $\gamma_x^-(t) = x + i \delta e^{it}$ for $0 \le t \le \pi$.
Calculating a residue we conclude
\begin{equation} \label{Residue}
\frac{f_{\alpha}(x)}{f_{\beta}(x) f_{\gamma}(x)} =
\frac{\mathcal{R}_\infty(x;\alpha,\beta,\xi)}{f_\xi(x)} +
\frac{\mathcal{R}_\infty(-x;-\alpha,-\xi, -\beta)}{f_\beta(x)}.
\end{equation}
The non-terminating $q$-Saalsch\"{u}tz sum follows easily from
(\ref{New}) and (\ref{Residue}) in the form
\begin{eqnarray*}
\frac{(\alpha u, \alpha/u;q)_\infty}{(\beta u, \beta/u, \xi u, \xi
/ u;q)_\infty} & = & \frac{(\alpha \xi,
\alpha/\xi;q)_\infty}{(\beta \xi, \beta/\xi, \xi u, \xi /
u;q)_\infty} \, \ _3\phi_2 \Big[
\begin{array}{c} \alpha/\beta, \ \xi u, \ \xi/u \\ \alpha \xi, \ q
\xi/ \beta \end{array}; \ q, \ q \Big] \\ & + & \frac{(\alpha
\beta, \alpha/\beta;q)_\infty}{(\beta u, \beta/u, \xi \beta, \xi /
\beta;q)_\infty} \ _3\phi_2 \Big[
\begin{array}{c} \alpha/\xi, \ \beta u, \ \beta/u \\ \alpha \beta, \ q
\beta / \xi \end{array}; \ q, \ q \Big].
\end{eqnarray*}
Note that our restrictions on $\alpha, \beta$ and $\xi$ disappear
after analytic extension.

\begin{remark}
Going back to the example considered in Paragraph
\ref{Subsection-Geometric}, we may consider the $\Po$-sequence
$z_t = q^{-t} \xi$ and apply Theorem
\ref{Theorem-Integral-Remainder} instead of Theorem
\ref{Teorema-Convergencia-Series}. In that case, similar arguments
give rise to a \emph{non-symmetrized} version of the
non-terminating $q$-Vandermonde sum $$\frac{(\alpha
x;q)_\infty}{(\beta x;q)_\infty} = \frac{(\alpha
\xi;q)_\infty}{(\beta \xi;q)_\infty} \ _2\phi_1 \Big[
\begin{array}{c} \alpha/\beta, \ x/\xi \\ q/\beta \xi, \end{array};
\ q, \ q \Big] + \frac{1}{2\pi i} \int_{-i \infty}^{i \infty}
\frac{(\alpha y, x /\xi;q)_\infty}{(\beta y, y / \xi;q)_\infty} \,
\frac{dy}{y-x}.$$
\end{remark}

\section{Binomial theorem}
\label{Section4}

As announced in the Introduction, we finish with a list of
binomial type formulas adapted to our classification into
canonical forms. We must recall that the same techniques we are
using here were already employed by Ismail in \cite{I} to obtain
the binomial type sum for the trigonometric canonical form, see
below. Let $f \in \Pol_n[x]$ and $(z_t)$ a $\Po$-sequence with
$z_0, z_1, \ldots z_n$ pairwise distinct. Then, it follows from
Remark \ref{Observacion-Polinomio-Interpolador} that
$\mathcal{R}_m f = 0$ for any $m \ge n$. Moreover, by a simple
continuity argument, we may drop the assumption that $z_0, z_1,
\ldots z_n$ are pairwise distinct. We need to use the $q$-binomial
coefficients $$\Big[ \!\!
\begin{array}{c} r
\\ k \end{array} \!\! \Big]_q = \frac{(q;q)_r}{(q;q)_k (q;q)_{r-k}}
\qquad \mbox{where} \qquad (x;q)_n = \prod_{k=0}^{n-1} (1 - q^k
x).$$

\begin{theorem} \label{Teorema-Binomial-Theorem}
Given two $\Po$-sequences $(y_t)$ and $(z_t)$, we have
$$\prod_{k=0}^{n-1} (x-z_k) = \sum_{k=0}^n \Big[ \!\!
\begin{array}{c} n \\ k \end{array} \!\! \Big]_{q} q^{-k (n-k)/2}
\prod_{j=0}^{k-1} (x-y_j) \prod_{j=0}^{n-k-1} (y_{k/2} - z_{j +
\frac{k}{2}}).$$
\end{theorem}

\begin{proof}
We claim that the following identity holds $$\partial_k
\Phi_n(x,\cdot) = \Big( \prod_{j=0}^{k-1} \frac{\lambda^{n-j} -
\lambda^{j-n}}{\lambda^{k-j} - \lambda^{j-k}} \Big)
\Phi_{n-k}(x,\cdot), \qquad 0 \le k \le n.$$ It can be checked by
induction on $k$ and applying the last identity in the proof of
Lemma \ref{Lema-Dk-Analitico}. In particular, since
$\Phi_n(z_{(n-1)/2},\cdot) \in \Pol_n[x]$, Theorem
\ref{Teorema-Taylor-Series} gives $$\Phi_n(z_{(n-1)/2},x) =
\sum_{k=0}^n \Big( \prod_{j=0}^{k-1} \frac{\lambda^{n-j} -
\lambda^{j-n}}{\lambda^{k-j} - \lambda^{j-k}} \Big)
\Phi_{n-k}(z_{(n-1)/2},y_{k/2}) \, \prod_{j=0}^{k-1} (x - y_j).$$
By the definition of $\Phi_k$ and $q$-binomial coefficient, the
proof is completed.
\end{proof}

Theorem \ref{Teorema-Binomial-Theorem} is symmetric under the
transformation $q \mapsto q^{-1}$. Moreover, when $q=1$ we must
take the obvious limits. Clearly, for the continuous canonical
form, Theorem \ref{Teorema-Binomial-Theorem} is nothing but
Newton's binomial theorem. Other explicit formulas of binomial
type arise from the remaining canonical forms. Moreover, them can
be rewritten in terms of the basic hypergeometric function. In the
following table we summarize the binomial type expressions and we
mention the corresponding hypergeometric formula.

\footnotesize
\begin{center}
\setlength{\extrarowheight}{7pt}
\begin{tabular}{|r|c|c|} \multicolumn{1}{c}{} &
\multicolumn{1}{c}{\normalsize \textbf{Binomial form}} &
\multicolumn{1}{c}{\normalsize\textbf{Sum}}
\\ \cline{1-3} \textbf{C} & $\displaystyle (x-z)^n = \sum_{k=0}^n
\binom{n}{k} (x-y)^k (y-z)^{n-k}$ & \textbf{\emph{Newton}}
\\ \cline{1-3} \textbf{A} & $\displaystyle \binom{x-z}{n} =
\sum_{k=0}^n \binom{x-y}{k} \binom{y-z}{n-k}$ &
\textbf{\emph{Gauss}$^{\dag}$} \\ \cline{1-3} \textbf{Q} &
$\displaystyle (u+w)_n \binom{u-w}{n} = \sum_{k=0}^n
\binom{u-v}{k} (u+v)_k \binom{v-w}{n-k} (v+w+k)_{n-k}$ &
\textbf{\emph{Pfaff}} \\ \cline{1-3} \textbf{G} & $\displaystyle
\prod_{k=0}^{n-1} (x - q^k z) = \sum_{k=0}^n {{\, n \,}\brack{\, k
\,}}_q \prod_{j=0}^{k-1} (x - q^j y) \prod_{j=0}^{n-k-1} (y - q^j
z)$ & \textbf{\emph{$q$-Gauss}$^{\dag}$} \\ \cline{1-3} \textbf{T}
& $\displaystyle \frac{(wu,w/u;q)_n}{(q;q)_n} = \sum_{k=0}^n
\frac{(vu,v/u;q)_k}{(q;q)_k} \frac{(wv
q^k,w/v;q)_{n-k}}{(q;q)_{n-k}} \Big( \frac{w}{v} \Big)^k$ &
\textbf{\emph{$q$-Saalsch\"{u}tz}}
\\ \cline{1-3} \end{tabular} \\ \vskip10pt \normalsize
\textsc{Table II.} Binomial type summation formulas. \vskip6pt
\end{center}

\normalsize

The reader is referred to Gasper and Rahman's book \cite{GR} to
see these summation formulas. The symbol $\dag$ means that, if we
rewrite the binomial type formula in terms of hypergeometric
summation formulas, what we obtain is a terminating form of the
corresponding non-terminating summation formula. Although we leave
the details for the interested reader, we should point out the
changes of variables employed in the quadratic and trigonometric
canonical forms.
$$\begin{array}{clll} \mathbf{(Q)} & x = u^2 & y_k = (v+k)^2 & z_k
= (w+k)^2 \\ \mathbf{(T)} & x = \frac{1}{2} (u + u^{-1}) & y_k =
\frac{1}{2} (q^k v + q^{-k} v^{-1}) & z_k = \frac{1}{2} (q^k w +
q^{-k} w^{-1}).
\end{array}$$ In the trigonometric case, we also need to use the
factorization formula
(\ref{Ecuacion-Factorizacion}).

\section*{Acknowledgments}

The authors wish to thank the referee for several comments that
gave rise to a significant improvement of this paper. This
research was supported in part by the Project MTM2004-00678,
Spain.

\enlargethispage{1cm}

\bibliographystyle{amsplain}

\end{document}